\documentclass[12pt]{article}
\usepackage{amsmath}
\usepackage{amsfonts}
\usepackage{amssymb}
\usepackage{amsthm}
\usepackage{hyperref}

\usepackage{graphicx}

\def\bee{\begin{equation}}
\def\eee{\end{equation}}

\begin{document}

\thispagestyle{empty}
\centerline{}
\bigskip
\bigskip
\bigskip
\bigskip
\bigskip
\bigskip
\centerline{\Large\bf Two  arguments that the nontrivial zeros}
\bigskip
\centerline{\Large\bf of the Riemann zeta function are irrational}
\bigskip
\bigskip

\begin{center}
{\large \sl Marek Wolf}\\*[5mm]
e-mail:mwolf@ift.uni.wroc.pl\\
\bigskip

\end{center}

\bigskip
\bigskip

\begin{center}
{\bf Abstract}\\
\bigskip
\begin{minipage}{12.8cm}
We have used the first 2600 nontrivial
zeros $\gamma_l$ of the Riemann zeta function calculated  with 1000 digits accuracy
and developed them  into
the continued fractions. We  calculated the geometrical means of the denominators of
these continued fractions 
and for all cases we get values close to the
Khinchin's constant, what suggests that $\gamma_l$ are irrational. Next we
have calculated the $n$-th square roots of the denominators $q_n$ of the
convergents of the continued fractions obtaining values close to the Khinchin–--L{\`e}vy
constant, again supporting the common believe that $\gamma_l$ are irrational.
\end{minipage}
\end{center}

\bigskip\bigskip

\bibliographystyle{abbrv}

\section{Introduction}

Bernhard G.F. Riemann has shown \cite{Riemann} that the l.h.s. of the identity
valid only for $\Re[s] >1 $:
\bee
\zeta(s) := \sum_{n=1}^\infty \frac{1}{n^s} = \prod_{p=2}^\infty \left(1-\frac{1}{p^s}\right)^{-1},
~~~~s=\sigma+it
\eee
can be analytically continued to the whole
complex plane devoid of $s=1$ by means of the following contour integral:
\bee
\zeta(s)=\frac{\Gamma(-s)}{2\pi i}\underset{\mathcal{P}}\int            
\frac{(-x)^s}{e^x-1} \frac{dx}{x}
\eee
\noindent where the integration is performed along the path ${\mathcal{P}}$

\vskip 0.3cm

\begin{picture} (0,120)(0,0) \thicklines
  \put(110,40){\vector(0,1){80}} \put(50,80){\vector(1,0){120}}

\put(110,80){\oval(40,40)[l]}
\put(110,85){\oval(30,30)[rt]}
\put(110,75){\oval(30,30)[rb]}

\put(125,75){\vector( 1,0){45}}
\put(170,85){\vector(-1,0){45}}
\end{picture}

\vskip -1.2cm
Now dozens of integrals 
and series representing the  $\zeta(s)$
function are known,  for collection of  such formulas see for example  the entry
{\it Riemann Zeta Function} in \cite{Weisstein} and references cited therein.

The $\zeta(s)$ function has trivial zeros $-2, -4, -6, \ldots$ and infinity of
nontrivial complex zeros $\rho_l=\beta_l+i\gamma_l$ in the critical strip: $\beta_l\in (0,1)$.
The Riemann Hypothesis (RH) asserts that $\beta_l=\frac{1}{2}$  for all $l$  --- i.e.
all zero lie on the critical line $\Re(s)=\frac{1}{2}$.  Presently it is added  that
these nontrivial zeros are simple: $\zeta'(\rho_l)\neq 0$ --- many explicit
formulas of number theory contain $\zeta'(\rho_l)$ in the denominators.
In  1914 G. Hardy \cite{Hardy1914} proved that infinitely many zeros of $\zeta(s)$ lie on
the critical line. A. Selberg  \cite{Selberg1942} in 1942 has shown that at least a
(small) positive proportion of the zeros of $\zeta(s)$ lie on the critical line.
The first quantitative  result was
obtained by  N. Levinson in 1974 \cite{Levinson1974} who showed that at least one-third
of the zeros lie on the critical line.  In 1989 B. Conrey  \cite{Conrey1989}
improved this to two-fifths and quite recently with collaborators \cite{Conrey2010}
to over 41\%.    It was checked computationally \cite{Gourdon}
that the $10^{13}$
first zeros of the Riemann Zeta function fulfill the condition $\beta_l=\frac{1}{2}$.
A. Odlyzko checked that RH is true in different intervals
around $10^{20}$ \cite{Odlyzko_zera1}, $10^{21}$ \cite{Odlyzko_zera2}, $10^{22}$
\cite{Odlyzko2001}, see also \cite{Gourdon} for the two billion zeros from the
zero $10^{24}$.

There is no hope to obtain the analytical formulas 
for the imaginary parts $\gamma_l$ of
the nontrivial zeros of $\zeta(s)$ but the common belief is that they are irrational
and perhaps even transcendental \cite{Odlyzko-chaos}.
The problem of any linear relations between $\gamma_l$  with integral
coefficients  appeared for the first time in the paper of  A.E. Ingham \cite{Ingham1942}
in connection with the  Mertens conjecture. This conjecture specifies the growth of the
function  $M(x)$  defined by
\bee
M(x)=\sum_{n\leq x} \mu(n), 	
\eee
where $\mu(n)$ is the M{\"o}bius function
\bee
\mu(n) \,=\,
\left\{
\begin{array}{ll}
1 & \mbox {for  $ n =1 $} \\
0 & \mbox {when $p^2|n$}\\
(-1)^r & \mbox{\rm when}~ n=p_1 p_2 \ldots p_r
\end{array}
\right.
\eee
The Mertens conjecture claims that
\bee
|M(x)|<x^\frac{1}{2}.
\label{Mertens}
\eee
From this inequality  the RH would follow.   A. E. Ingham  in  \cite{Ingham1942} showed
that the validity of the Merten's conjecture requires that the imaginary parts of
the nontrivial zeros should fulfill the relations of the form:
\bee
\sum_{l=1}^N c_l \gamma_l=0,
\eee
where $c_l$ are integers not all equal to zero. This result raised the doubts in the
inequality (\ref{Mertens})  and indeed in 1985 A. Odlyzko and H. te Riele
\cite{odlyzko1985} disproved  the Merten's conjecture.

In this paper we are going to exploit two facts about the continued fractions: the
existence of the   Khinchin constant and  Khinchin–--L{\`e}vy constant,
see e.g. \cite[\S 1.8]{Finch},  to support the irrationality of $\gamma_l$.  Let
\bee
r = [a_0(r); a_1(r), a_2(r), a_3(r), \ldots]=
a_0(r)+\cfrac{1}{a_1(r) + \cfrac{1}{a_2(r) + \cfrac{1}{a_3(r)+\ddots}}}
\eee
be the continued fraction expansion of the real number $r$, where $a_0(r)$ is an integer and all $a_k(r)$ with $k\geq 1$  are positive integers.
Khinchin has proved \cite{Khinchin}, see also \cite{Ryll-Nardzewski1951}, that
\bee
\lim_{n\rightarrow \infty} \big(a_1(r) \ldots a_n(r)\big)^{\frac{1}{n}}=
\prod_{m=1}^\infty {\left\{ 1+{1\over m(m+2)}\right\}}^{\log_2 m} \equiv K_0 \approx 2.685452001\dots
\label{Khinchin}
\eee
is a constant for almost all real $r$ \cite[\S 1.8]{Finch}. The exceptions are of the
Lebesgue measure zero and include {\it rational numbers}, quadratic irrationals and
some irrational numbers too, like for example the
Euler constant $e=2.7182818285\ldots$ for which the limit (\ref{Khinchin}) is
infinity.  The constant $K_0$ is called the Khinchin constant.
If the quantities
\bee
K(r; n)=\big(a_1(r) a_2(r) \ldots a_n(r)\big)^{\frac{1}{n}}
\label{Kny}
\eee
for a given number $r$ are close to $K_0$ we can regard it as an indication
that $r$ is irrational.

Let the rational $p_n/q_n$ be the $n$-th  partial convergent of the continued fraction:
\bee
\frac{p_n}{q_n}=[a_0; a_1, a_2, a_3, \ldots, a_n].
\eee
For almost all real numbers $r$ the denominators of the  finite continued fraction
approximations fulfill:
\bee
\lim_{n \rightarrow \infty} \big(q_n(r)\big)^{1/n} = e^{\pi^2/12\ln2} \equiv L_0 = 3.275822918721811\ldots
\eee
where $L_0$ is called the  Khinchin---L{\`e}vy's constant  \cite[\S 1.8]{Finch}.
Again the set of exceptions to
the above limit is  of the Lebesgue measure zero and it includes rational numbers,
quadratic irrational etc.

\section{The computer experiments}

First 100 zeros $\gamma_l$  of $\zeta(s)$  accurate to over 1000 decimal places we have taken from
\cite{Odlyzko_zeros}.
Next 2500 zeros of $\zeta(s)$   with precision of 1000 digits
were  calculated using the built in Mathematica v.7 procedure \verb"ZetaZero[m]".
We have checked  using PARI/GP \cite{PARI} that these zeros were accurate within at
least 996
places in the sense that in the worst case $|\zeta(\rho_l)|<10^{-996}, ~l= 1, 2, \ldots, 2600$.
PARI has built in function
\verb"contfrac"$(r,\{nmax\})$ which creates the row vector ${\bf a}(r)$ whose
components are the denominators $a_n(r)$  of the continued fraction
expansion of $r$, i.e.  ${\bf a}=[a_0(r); a_1(r), \dots,a_n(r)]$ means that
\bee
r \approx
a_0(r)+\cfrac{1}{a_1(r) + \cfrac{1}{a_2(r) + \cfrac{1}{\ddots\cfrac{1}{a_n(r)}}}}
\eee
The parameter $nmax$ limits the number of terms $a_{nmax}(r)$; if it is omitted
the expansion stops with a declared precision
of representation of real numbers at the last significant partial quotient.

By trials we have determined that the precision set to {\tt $\backslash$p  2200} is sufficient in the
sense that scripts with larger precision generated exactly the same results: the rows
${\bf a}(\gamma_l)$ obtained  with accuracy 2200 digits were the same for
all $l$  as those
obtained for accuracy 2600 and the continued fractions  accuracy set to 2100 digits
gave different denominators $a_n(\gamma_l)$  
With the precision set  to 2200 digits we have developed the 1000 digits values of
 each $\gamma_l$, $l=1, 2, \ldots 2600$, into the continued fractions
\bee
\gamma_l \doteq [a_0(l); a_1(l), a_2(l), a_3(l), \ldots, a_{n(l)}(l)]\equiv {\bf a}(l)
\label{g-cfr}
\eee
without specifying the
parameter $nmax$, thus the length of the vector ${\bf a}(l)$  depended on $\gamma_l$
and it turns out that the number of denominators was contained between 1788 and  2072.
The value of the product  $a_1 a_2 \ldots a_{n(l)} $  was typically of the
order $10^{800}-10^{870}$.
Next for each $l$ we have calculated the geometrical means:
\bee
K_l(n(l))= \left( \prod_{k=1}^{n(l)} a_k(l) \right)^{1/n(l)}.
\eee
The results are presented in the Fig.1. Values of $K_l(n(l))$ are scattered around the
red line representing $K_0$.  To gain some insight into the rate of convergence
of $K_l(n(l))$ we have plotted in the Fig. 2 the number of sign changes $S(l)$ of
$K_l(m)-K_0$ for each $l$ when $m=100, 101, \ldots n(l)$, i.e.
\bee
S_K(l)={\rm number ~ of ~ such ~{\it  m}~ that} ~~~(K_l(m+1)-K_0)(K_l(m)-K_0)<0.
\eee
The largest $S_K(l)$ was  122 and it occurred for the zero $\gamma_{194}$  and for 381
zeros there were no sign changes at all. In the Fig. 3 we present plots of
$K_l(m)$ as a function of $m$  for a few zeros $\gamma_l$.

Let the rational $p_{n(l)}(\gamma_l)/q_{n(l)}(\gamma_l)$ be the $n$-th
partial convergent of the continued fractions (\ref{g-cfr}):
\bee
\frac{p_{n(l)}(\gamma_l)}{q_{n(l)}(\gamma_l)}={\bf a}(l).
\eee
For each zero $\gamma_l$ we have calculated the partial convergents
$p_{n(l)}(\gamma_l)/q_{n(l)}(\gamma_l)$.
Next from these denominators $q_{n(l)}(\gamma_l)$ we have calculated the quantities
$L_l(n(l))$:
\bee
L_l(n(l))= \left(q_{n(l)}\right)^{1/n(l)}, ~~~l=1, 2, \ldots , 2600
\eee
The obtained values of $L_l(n(l))$  are presented in the Fig.4.
These  values scatter around the red line representing the Khinchin---L{\`e}vy's
constant $L_0$. As in the case of $K_l(m)$ the Fig.5 presents the number of
sign changes of the difference $L_l(m)-L_0$ as a function of the index $m$
of the denominator of the $m$-th convergent $p_m/q_m$
\bee
S_L(l)={\rm number ~ of ~ such ~ {\it m}~ that} ~~~(L_l(m+1)-L_0)(L_l(m)-L_0)<0.
\eee
The maximal number of sign
changes was 136 and appeared for the zero $\gamma_{1389}$ and there were 396 zeros
without sign changes.

In the Fig. 6 we have plotted the ``running'' absolute difference between $K_l(m)$ and
$K_0$  averaged over all 2600 zeros:
\bee
A_K(m)=\frac{1}{2600}\sum_{l=1}^{2600} \big | K_l(m)-K_0\big|,~~m=100,101, \ldots 1788
\eee
and the similar  average for the difference between $L_l(m)$ and $L_0$:
\bee
A_L(m)=\frac{1}{2600}\sum_{l=1}^{2600} \big |L_l(m)-L_0\big|,~~m=100,101, \ldots 1788.
\eee
These two averages very rapidly tend to zero. Although it does not prove nothing,
the fact that  the curves representing $A_K(m)$ and $A_L(m)$ almost coincide
is very convincing.  In the inset the plot on double logarithmic scale reveals that
both $A_L(m)$ and $A_k(m)$ decrease like $C_{K,L}/\sqrt{m}$ where $C_K= 2.3868\ldots$
and $C_L=2.4473\ldots$.  It is a pure speculation  linking  the power of these
dependencies $m^{-1/2}$ to the ordinate of the critical line.

\section{Concluding remarks}

There are generalizations of above quantities $K(n)$ given by (\ref{Kny}).
It can be shown that the following $s$-mean values of the denominators $a_k(r)$
of the continued fraction for a real number $r$:
\bee
M(n,s;r)= \left(\frac{1}{n} \sum_{k=1}^n \big( a_k(r)\big)^s \right)^{1/s}
\label{general-M}
\eee
are divergent for $s\geq 1$ and convergent for $s<1$ for almost all real $r$
\cite[\S 1.8]{Finch}.  It can be shown that for $s<1$
\bee
\lim_{n\to\infty} M(n, s;r)=\left(\sum_{k=1}^\infty -k^s \log_2\left( 1-\frac{1}{(k+1)^2} \right) \right)^{1/s} \equiv K_s
\eee
where $M(n, s;r)$ for almost all $r$ are the same.
The quantities (\ref{general-M}) can be computed for imaginary parts
of nontrivial zeta zeros $M(n,s; \gamma_l)$ and compared with values of $K_s$
but we leave it for further investigation.

The continued fraction were used in the past in the  Apery`s  proof of
irrationality of $\zeta(3)$.  In the paper  \cite{Cvijovic-et-al-1997}
values  of $\zeta(n)$ for all $n\geq 2$ were expressed in terms of rapidly
converging continued fractions. These results were analytical, but in case
of the nontrivial zeros of the $\zeta(s)$ function  we are left only
with the computer experiments. The results reported in this paper suggests
that they are irrational.

\begin{figure}
\vspace{-0.3cm}
\begin{center}
\includegraphics[width=1.2\textwidth,angle=90]{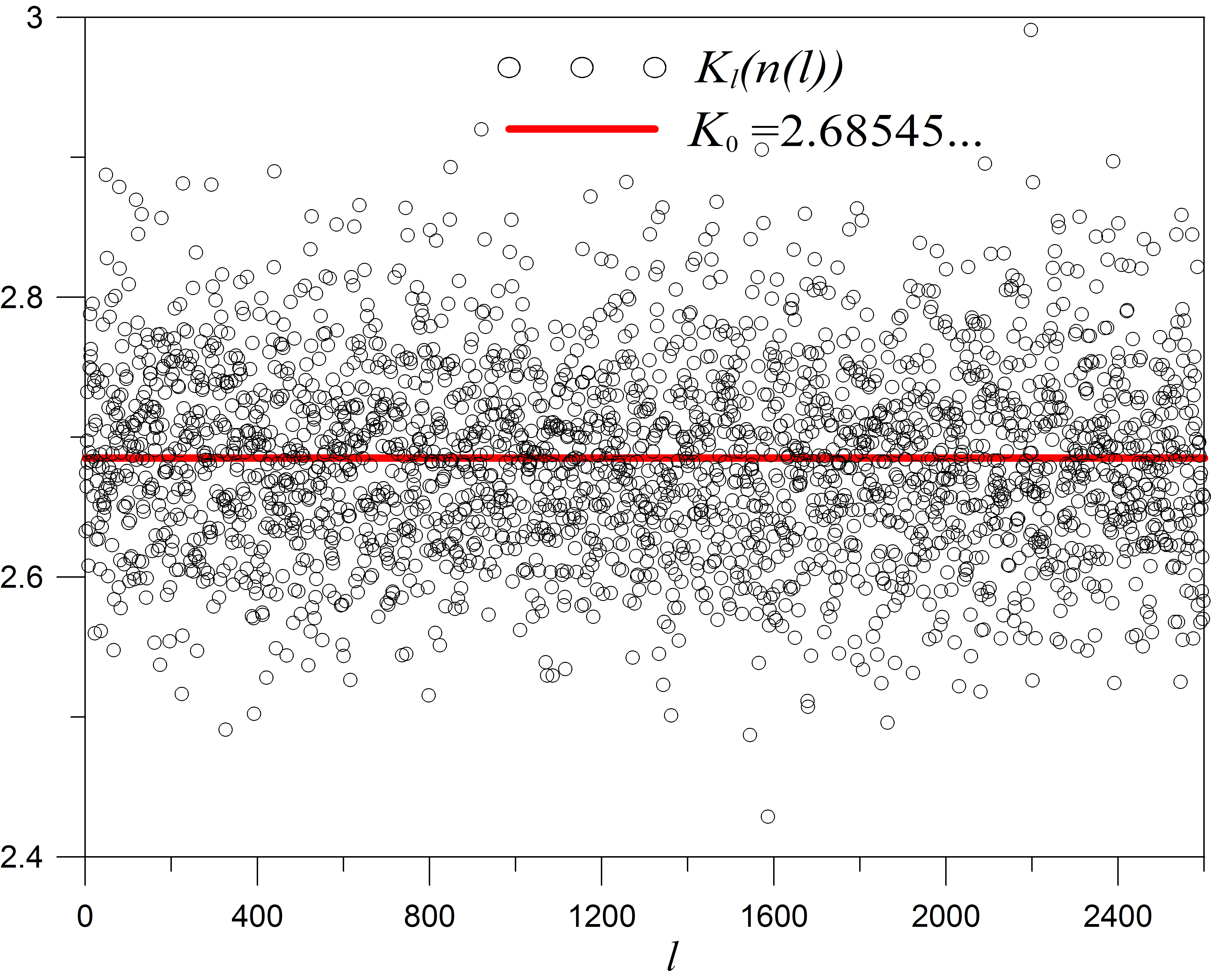} \\
\vspace{1.7cm} Fig.1 The plot of $K_l(n(l))$. \\
\end{center}
\end{figure}

\begin{figure}
\vspace{-0.3cm}
\begin{center}
\includegraphics[width=1.2\textwidth,angle=90]{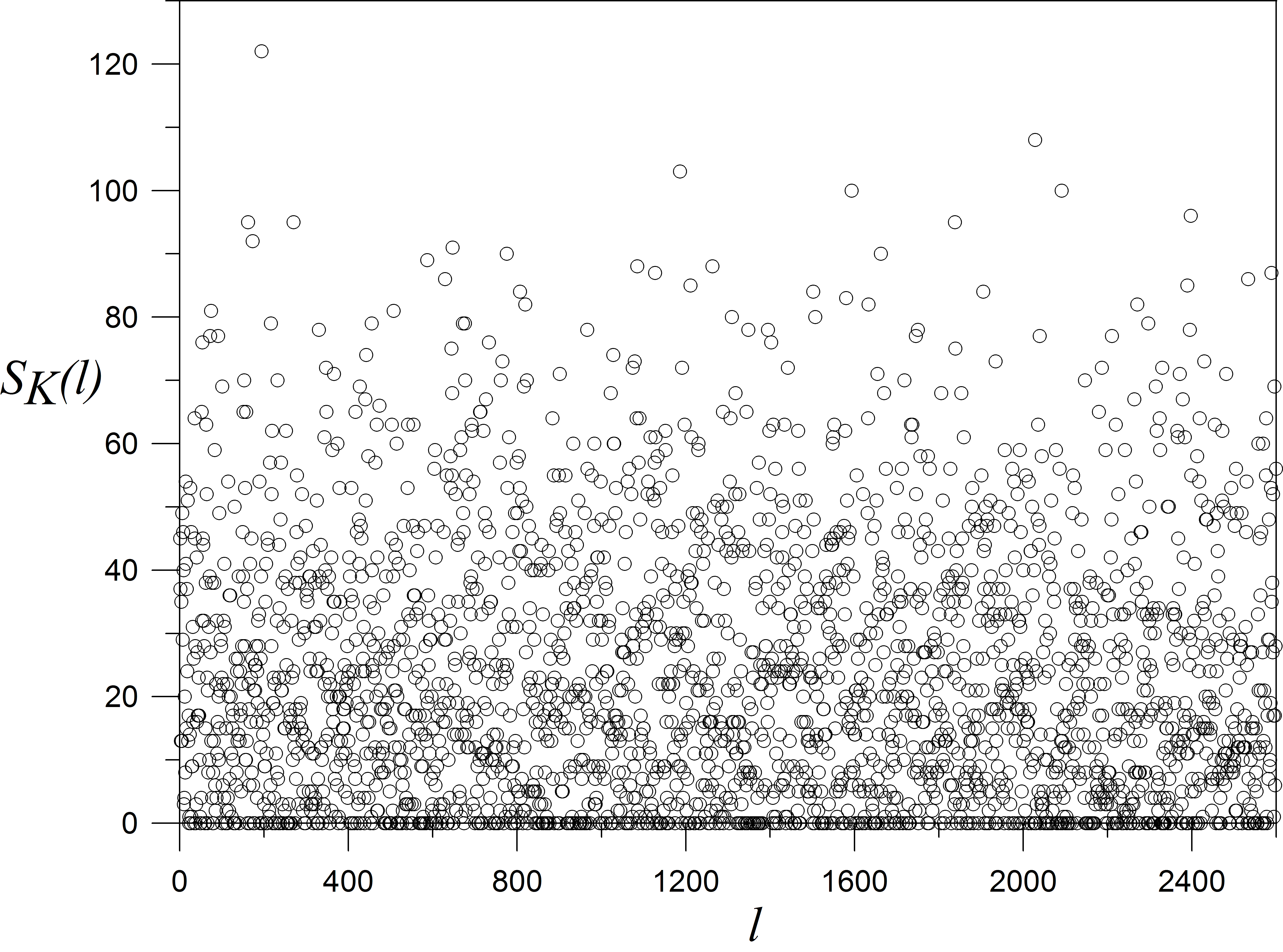} \\
\vspace{1.7cm} Fig.2 The number of  such $m$ that 
$(K_l(m+1)-K_0)(K_l(m)-K_0)<0$ for each $l$ (the initial transient values
of $m$ were skipped--- --- sign changes were detected for $m=100, 101, \ldots n(l)$).  \\
\end{center}
\end{figure}

\begin{figure}
\vspace{-0.3cm}
\begin{center}
\includegraphics[width=1.2\textwidth,angle=90]{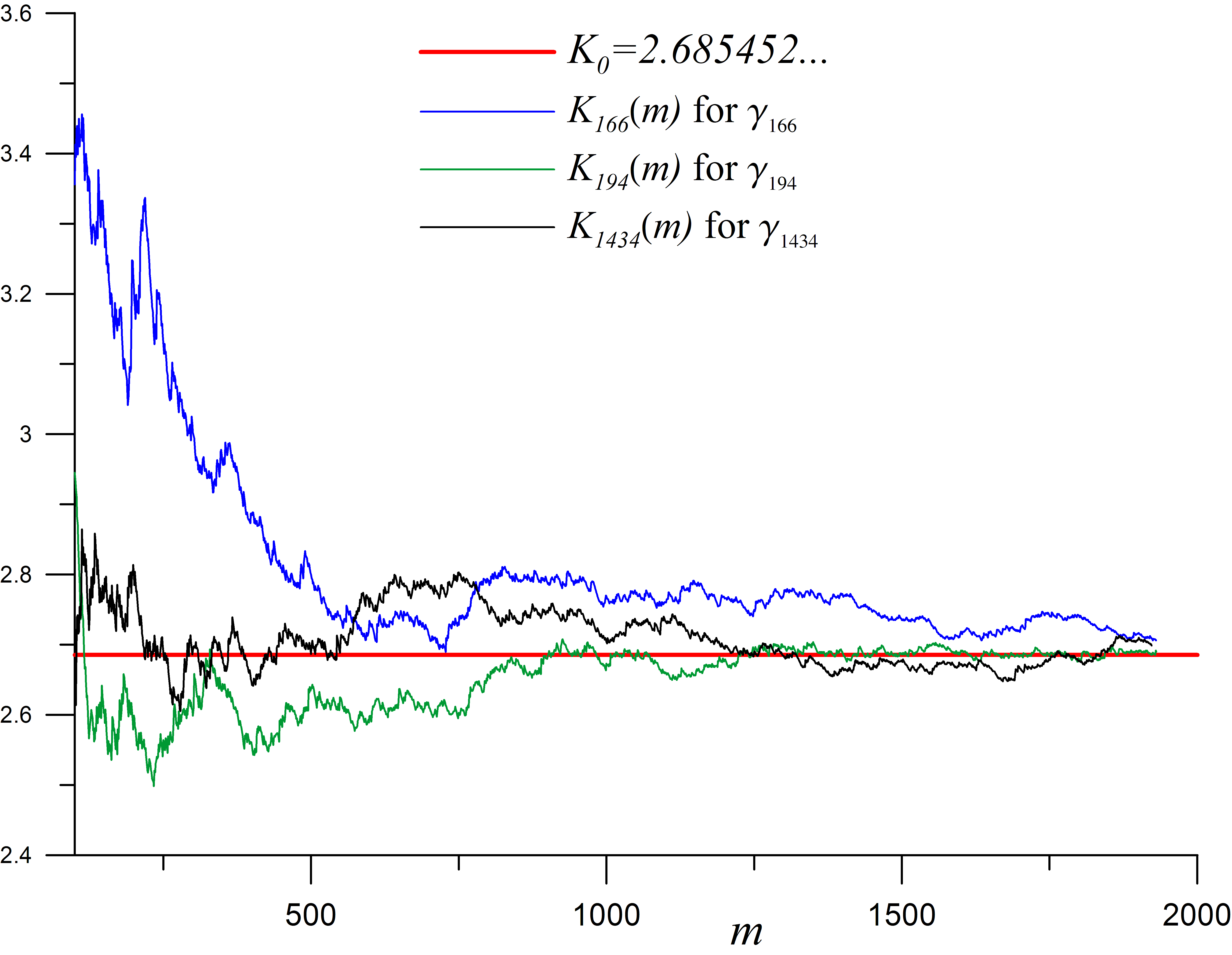} \\
\vspace{1.7cm}Fig.3  For $\gamma_{166}$ there was no sign change of the difference $K_{166}(m)-K_0$.
 For $\gamma_{194}$ there were 122 sign changes of the difference $K_{194}(m)-K_0$ ---
it was  the largest number of sign changes among all zeros.
For $\gamma_{1434}$ there were 63  sign changes of the difference $K_{1434}(m)-K_0$. \\
\end{center}
\end{figure}

\begin{figure}
\vspace{-0.3cm}
\begin{center}
\includegraphics[width=1.2\textwidth,angle=90]{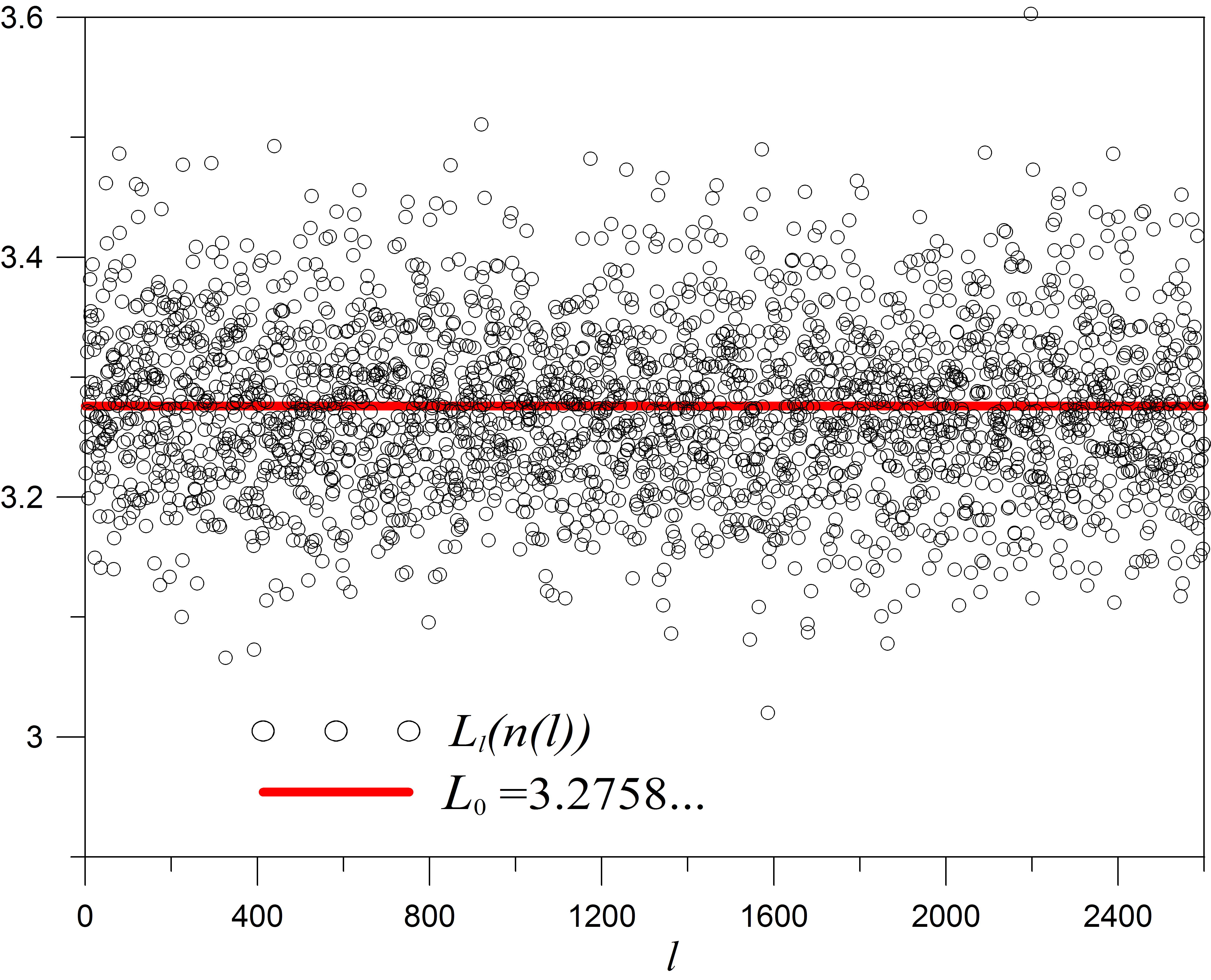} \\
\vspace{1.7cm} Fig.4 The plot of $L_l(n(l))$. \\
\end{center}
\end{figure}

\begin{figure}
\vspace{-0.3cm}
\begin{center}
\includegraphics[width=1.2\textwidth,angle=90]{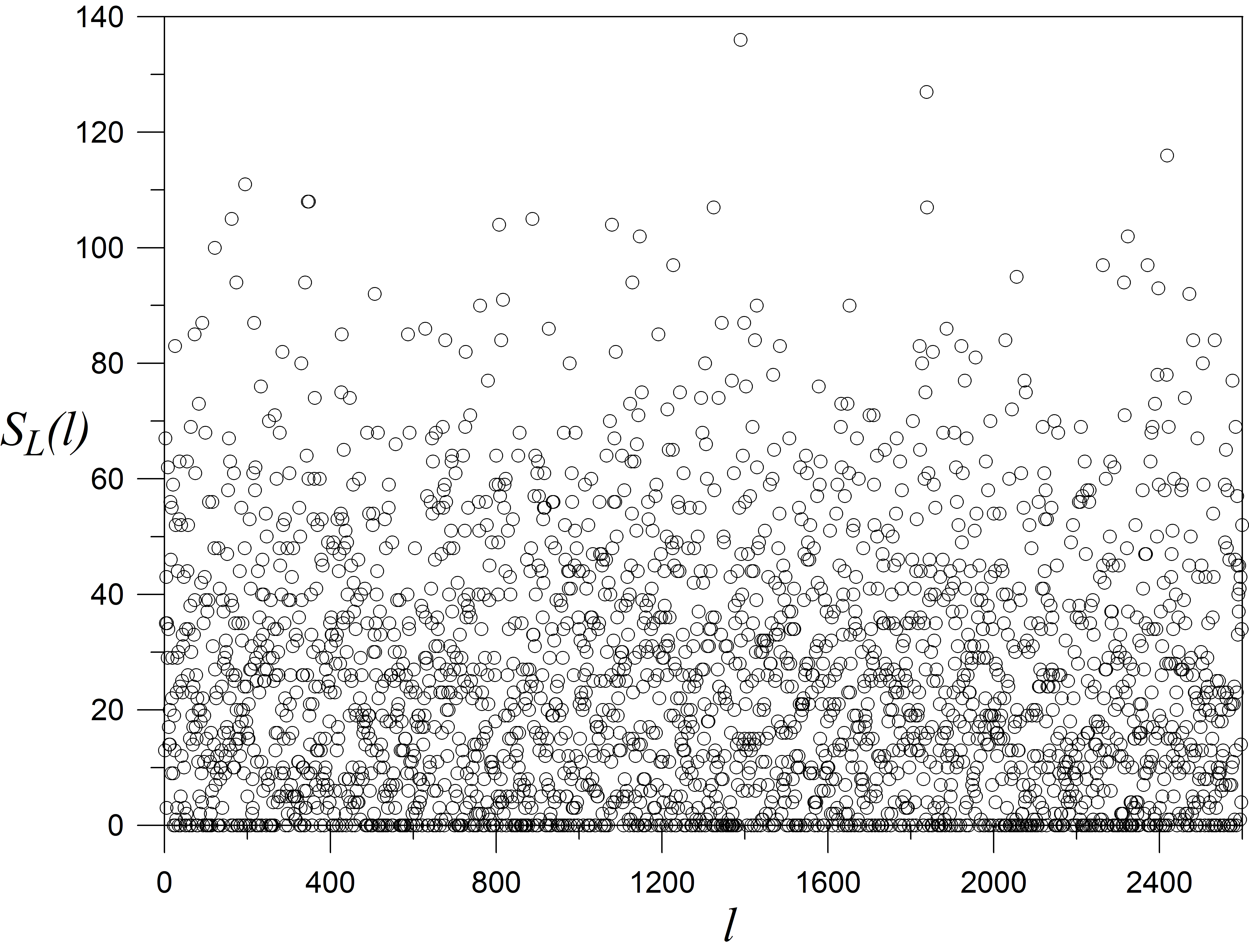} \\
\vspace{1.7cm} Fig.5 The number of  such $m$ that 
$(L_l(m+1)-L_0)(L_l(m)-L_0)<0$ for each $l$ (the initial transient values
of $m$ were skipped --- sign changes were detected for $m=100, 101, \ldots n(l)$).  \\
\end{center}
\end{figure}

\begin{figure}
\vspace{-0.3cm}
\begin{center}
\includegraphics[width=1.2\textwidth,angle=90]{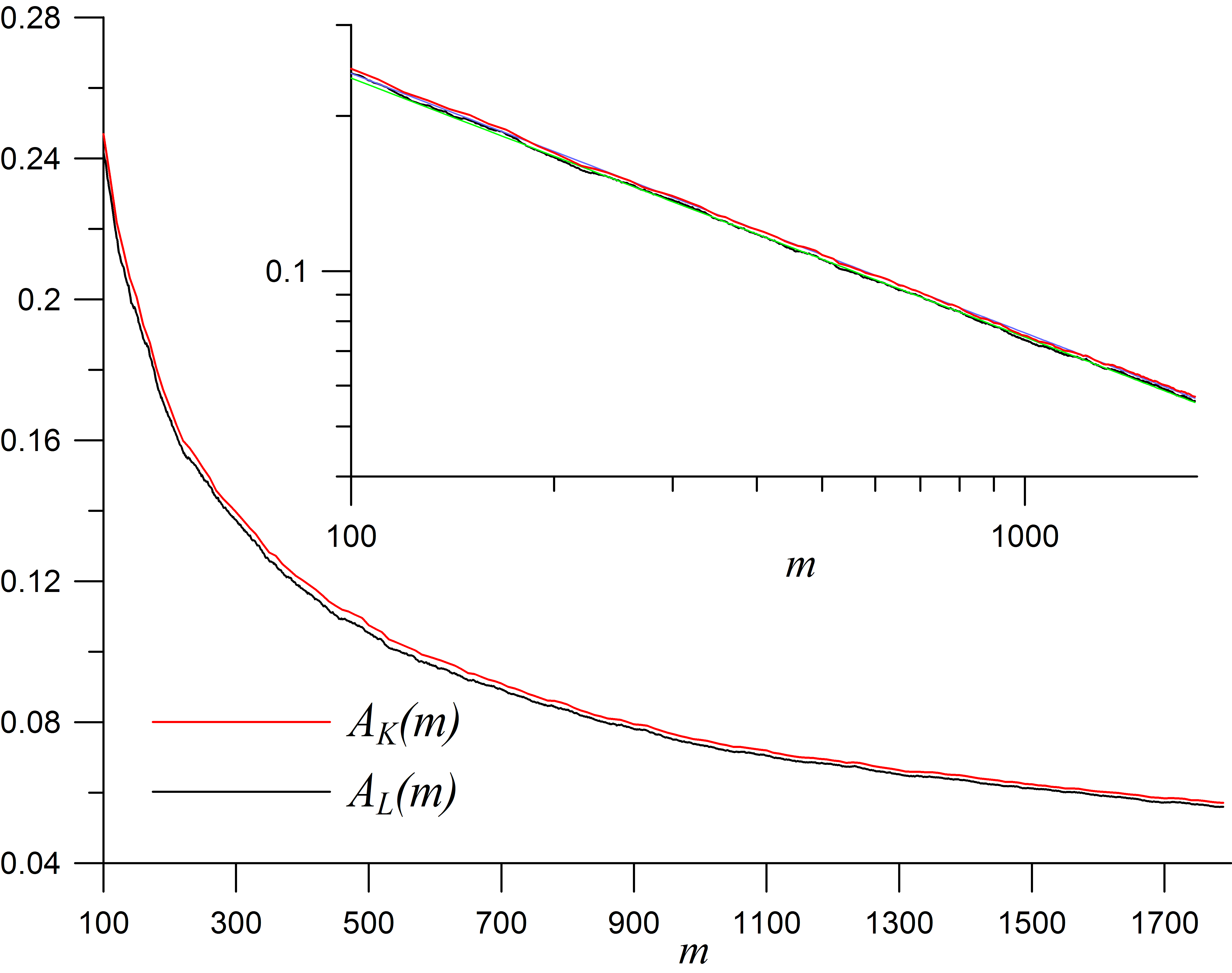} \\
\vspace{1.7cm} Fig.6 The averaged over all 2600 zeros differences $|K_l(m)-K_0|$
and  $|L_l(m)-L_0|$  plotted  for $m=100, 101, \ldots 1788$).  In the inset the same
curves are plotted on the double logarithmic scale together with fits obtained by the
least square method. The equations of the fits are $2.3868/m^{0.5019}$ for $A_K(m)$
(green line)  and $2.4473/m^{0.5028}$ for $ A_L(m)$ (blue line)\\
\end{center}
\end{figure}

\end{document}